\documentclass[10pt]{article}

\usepackage{amsmath}
\usepackage{amsthm}

\usepackage[bookmarks]{hyperref}
\hypersetup{
        colorlinks=true,
        linkcolor=black,
        anchorcolor=black,
        citecolor=black,
        urlcolor=black,
        pdfpagemode=UseThumbs,
        pdftitle={A probabilistic proof of a lemma that is not Burnside's},
        pdfauthor={Vincent Vatter},
}

\makeatletter

\newcommand\blfootnote[1]{%
  \begingroup
  \renewcommand\thefootnote{}\footnote{#1}%
  \addtocounter{footnote}{-1}%
  \endgroup
}

\let\saveqed\qed
\renewcommand\qed{%
   \ifmmode\displaymath@qed
   \else\saveqed
   \fi}

\title{\vspace{-1in}\sc A Probabilistic Proof of a\\ Lemma that is not Burnside's}
\author{%
Vincent Vatter\\[-0.25ex]
\small Department of Mathematics\\[-0.5ex]
\small University of Florida\\[-0.5ex]
\small Gainesville, Florida USA\\[-1.5ex]
}

\date{}

\newcommand{\st}{\::\:}
\newcommand{\prob}[1]{\operatorname{Pr}\bigl[\,#1\,\bigr]}
\newcommand{\fix}{\operatorname{fix}}

\newcommand{\orb}{\operatorname{orb}}
\newcommand{\probgiven}{\:\big\vert\:}

\begin{document}
\maketitle
\thispagestyle{empty}
\blfootnote{This note appeared as a ``MathBit'' in \emph{Amer. Math. Monthly} {\bf 127}(1) (2020), 63.}

\noindent If the group $G$ acts on the set $X$, we define the \emph{orbit} of the element $x\in X$ as
\[
	\orb(x)=\{gx\st g\in G\}.
\]
The orbits partition $X$, and the orbit-counting lemma (see Neumann~\cite{neumann:a-lemma-that-is:} for the history of its name) shows how to compute the number of orbits. Bogart~\cite{bogart:an-obvious-proo:} gave a proof of this result using only multisets and the product rule. We give a probabilistic version of his proof.

\newtheorem*{orbit-counting-lemma}{\rm\bf The Orbit-Counting Lemma}
\begin{orbit-counting-lemma}
Suppose the finite group $G$ acts on the finite set $X$ and let $\fix(g)=\{x\in X\st gx=x\}$. Then the number of orbits of $X$ is
\[
	\frac{1}{|G|}\sum_{g\in G} |\fix(g)|.
\]
\end{orbit-counting-lemma}
\begin{proof}
Choose, in order and each uniformly at random, an element $g\in G$, an orbit of $X$, and an element $y$ of this orbit. Because $gy=x$ for some $x\in X$,
\[
	1
	=
	\sum_{x\in X}\prob{gy=x}
	=
	\sum_{x\in X}\prob{y\in\orb(x)}
	\cdot
	\prob{gy=x\probgiven y\in\orb(x)}.
\]
If $y\in\orb(x)$, then there is some $k\in G$ such that $kx=y$, so the mapping $h\mapsto k^{-1}h^{-1}$ is a bijection between $\{h\in G\st hy=x\}$ and $\{h\in G\st hx=x\}$. By applying this (measure-preserving) mapping, it follows that, for every $x\in X$,
\[
	\prob{gy=x\probgiven y\in\orb(x)}
	=
	\prob{gx=x\probgiven y\in\orb(x)}
	=
	\prob{gx=x}.
\]
Also, since the orbit $y$ was chosen from was itself chosen uniformly at random,
\[
	\frac{1}{\#\,\text{orbits}}\sum_{x\in X} \prob{gx=x}
	=
	1.
\]
This implies that the number of orbits is equal to
\[
	\sum_{x\in X} \prob{gx=x}
	=
	\frac{1}{|G|}\bigl|\bigl\{(g,x)\in G\times X\st gx=x\bigr\}\bigr|
	=
	\frac{1}{|G|}\sum_{g\in G} |\fix(g)|,
\]
as claimed.
\end{proof}

\end{document}